\documentclass{amsart}
\usepackage{m-pictex}
\title[The Two-Phase Membrane Problem]{The Two-Phase Membrane Problem --
an Intersection-Comparison Approach to the Regularity at Branch Points}
\author[H. Shahgholian ]{Henrik Shahgholian}
\address{Department of Mathematics, Royal Institute of Technology,
100~44  Stockholm, Sweden}
\email{henriksh@math.kth.se}
\author[G.S. Weiss]{Georg S. Weiss}
\address{Graduate School of Mathematical Sciences,
University of Tokyo, 3-8-1 Komaba, Meguro-ku, Tokyo-to, 153-8914 Japan,}
\email{gw@ms.u-tokyo.ac.jp,}
\thanks{$2000$ {\it Mathematics Subject Classification.\/} Primary
35R35, Secondary 35J60.}
\thanks{H. Shahgholian has been partially supported  by the Swedish Research Council.
G.S. Weiss has been partially supported by a Grant-in-Aid for Scientific Research,
Ministry of Education, Japan}
\thanks{G.S.Weiss wishes to thank the G\"oran Gustafsson Foundation for visiting appointments to the Royal Inst. of Technology in Stockholm. }
\date{}
\theoremstyle{plain}
\newtheorem{theorem}{Theorem}[section]
\newtheorem{lemma}[theorem]{Lemma}
\newtheorem{proposition}[theorem]{Proposition}
 \theoremstyle{definition}
 \theoremstyle{example}

\theoremstyle{definition}

\numberwithin{equation}{section}
\def\R{{\bf R}}

\def\dist{\hbox{\rm dist}}

\begin{document}
\begin{abstract}
For the two-phase membrane problem
$ \Delta u = {\lambda_+\over 2} \chi_{\{u>0\}} \> - \>
{\lambda_-\over 2} \chi_{\{u<0\}}\> ,$ where $\lambda_+> 0$ and
$\lambda_->0\> ,$ we prove in two dimensions that
the free boundary is in a neighborhood of each ``branch point''
the union of two $C^1$-graphs. We also obtain a stability
result with respect to perturbations of the boundary data.
Our analysis uses an intersection-comparison approach
based on the Aleksandrov reflection.\\
In higher dimensions
we show that the free boundary has finite $(n-1)$-dimensional Hausdorff
measure.
\end{abstract}
\maketitle
\section{Introduction}

In this paper  we study the regularity of the
obstacle-problem-like equation
\begin{equation}\label{obst}
\Delta u = {\lambda_+\over 2} \chi_{\{u>0\}} \> - \>
{\lambda_-\over 2} \chi_{\{u<0\}}\;    \qquad \hbox{in } \Omega ,
\end{equation}
where $\lambda_+> 0, \lambda_->0$ and $\Omega \subset \R^n$ is a
given domain. Physically the equation arises for example as the
``two-phase membrane problem'': consider an elastic membrane
touching the phase boundary between two liquid/gaseous phases with
different viscosity, for example a water surface. If the membrane
is pulled away from the phase boundary in both phases, then the
equilibrium state can be described by equation (\ref{obst}).
\begin{figure}
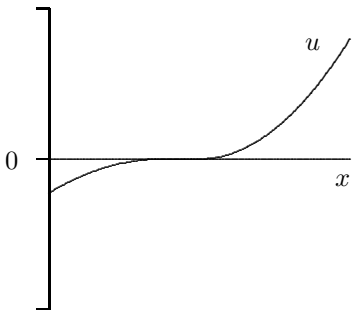

\[
\beginpicture
\setcoordinatesystem units <1cm,1cm>
\setplotarea x from -2 to 2, y from -2 to 2
\axis left ticks quantity 3 /
\put {\makebox(0,0)[t]{$u$}} [Bl] at 1.5 1.6
\put {\makebox(0,0)[t]{$x$}} [Bl] at 1.9 -0.2
\put {\makebox(0,0)[t]{$0$}} [Bl] at -2.5 0.1
\linethickness2pt
\plot -2 0.0 2 0.0 /
\linethickness20pt
\setquadratic
\plot -2 -0.45
-1.92 -0.40328
-1.84 -0.35912
-1.76 -0.31752
-1.68 -0.27848
-1.6 -0.242
-1.52 -0.20808
-1.44 -0.17672
-1.36 -0.14792
-1.28 -0.12168
-1.2 -0.098
-1.12 -0.07688
-1.04 -0.05832
-0.96 -0.04232
-0.88 -0.02888
-0.8 -0.018
-0.72 -0.00968
-0.64 -0.00392
-0.56 -0.00072
-0.48 0
-0.4 0
-0.32 0
-0.24 0
-0.16 0
-0.08 0
0 0
0.08 0.00256
0.16 0.01024
0.24 0.02304
0.32 0.04096
0.4 0.064
0.48 0.09216
0.56 0.12544
0.64 0.16384
0.72 0.20736
0.8 0.256
0.88 0.30976
0.96 0.36864
1.04 0.43264
1.12 0.50176
1.2 0.576
1.28 0.65536
1.36 0.73984
1.44 0.82944
1.52 0.92416
1.6 1.024
1.68 1.12896
1.76 1.23904
1.84 1.35424
1.92 1.47456
2 1.6 /
\endpicture
\]
\caption{A solution in 1d}\label{membrane}
\end{figure}

\begin{figure}
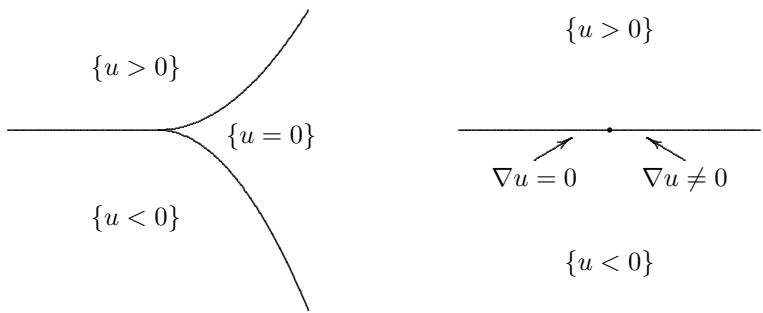

\[
\beginpicture
\setcoordinatesystem units <1cm,1cm>
\setplotarea x from -2 to 8, y from -2 to 2
\put {\makebox(0,0)[t]{$\{ u>0\}$}} [Bl] at -0.3 1
\put {\makebox(0,0)[t]{$\{ u=0\}$}} [Bl] at 1.5 0.1
\put {\makebox(0,0)[t]{$\{ u<0\}$}} [Bl] at -0.3 -1

\put {\makebox(0,0)[t]{$\{ u>0\}$}} [Bl] at 6 1.5
\put {\makebox(0,0)[t]{$\{ u<0\}$}} [Bl] at 6 -1.6
\put {\circle*{2}} [Bl] at 6 0
\put {\makebox(0,0)[t]{$\nabla u=0$}} [Bl] at 5 -.5
\put {\makebox(0,0)[t]{$\nabla u\ne0$}} [Bl] at 7 -.5
\arrow <2mm> [0.1,0.5] from 5 -.4 to 5.5 -0.1
\arrow <2mm> [0.1,0.5] from 7 -.4 to 6.5 -0.1

\linethickness2pt
\plot -2 0.0 0.0 0.0 /
\plot 4 0.0 8 0.0 /
\linethickness20pt
\setquadratic
\plot 0 0
0.04 0.00064
0.08 0.00256
0.12 0.00576
0.16 0.01024
0.2 0.016
0.24 0.02304
0.28 0.03136
0.32 0.04096
0.36 0.05184
0.4 0.064
0.44 0.07744
0.48 0.09216
0.52 0.10816
0.56 0.12544
0.6 0.144
0.64 0.16384
0.68 0.18496
0.72 0.20736
0.76 0.23104
0.8 0.256
0.84 0.28224
0.88 0.30976
0.92 0.33856
0.96 0.36864
1 0.4
1.04 0.43264
1.08 0.46656
1.12 0.50176
1.16 0.53824
1.2 0.576
1.24 0.61504
1.28 0.65536
1.32 0.69696
1.36 0.73984
1.4 0.784
1.44 0.82944
1.48 0.87616
1.52 0.92416
1.56 0.97344
1.6 1.024
1.64 1.07584
1.68 1.12896
1.72 1.18336
1.76 1.23904
1.8 1.296
1.84 1.35424
1.88 1.41376
1.92 1.47456
1.96 1.53664
2 1.6 /
\plot 0 0
0.04 -0.00096
0.08 -0.00384
0.12 -0.00864
0.16 -0.01536
0.2 -0.024
0.24 -0.03456
0.28 -0.04704
0.32 -0.06144
0.36 -0.07776
0.4 -0.096
0.44 -0.11616
0.48 -0.13824
0.52 -0.16224
0.56 -0.18816
0.6 -0.216
0.64 -0.24576
0.68 -0.27744
0.72 -0.31104
0.76 -0.34656
0.8 -0.384
0.84 -0.42336
0.88 -0.46464
0.92 -0.50784
0.96 -0.55296
1 -0.6
1.04 -0.64896
1.08 -0.69984
1.12 -0.75264
1.16 -0.80736
1.2 -0.864
1.24 -0.92256
1.28 -0.98304
1.32 -1.04544
1.36 -1.10976
1.4 -1.176
1.44 -1.24416
1.48 -1.31424
1.52 -1.38624
1.56 -1.46016
1.6 -1.536
1.64 -1.61376
1.68 -1.69344
1.72 -1.77504
1.76 -1.85856
1.8 -1.944
1.84 -2.03136
1.88 -2.12064
1.92 -2.21184
1.96 -2.30496
2 -2.4 /
\endpicture
\]
\caption{Examples of branch points}\label{branchp}
\end{figure}

Properties of the solution etc. have been derived by the authors in
\cite{interf} and \cite{uraltseva}. Moreover,
in \cite{SUW}, the authors gave a complete characterization of
global two-phase solutions
satisfying a quadratic growth condition at the two-phase free boundary point
$0$ and at
infinity. It turned out that each such solution coincides after
rotation with the one-dimensional solution $u(x) =
 {\lambda_+\over 4} \max(x_n,0)^2\> -{\lambda_-\over 4} \min(x_n,0)^2 .$
In particular this implies that each blow-up limit $u_0$ at
so-called ``branch points'', $\Omega \cap \partial\{ u>0\}\cap
\partial\{ u<0\}\cap \{\nabla u =0\}\> ,$ is after rotation of the
form $u_0(x) = {\lambda_+\over 4} \max(x_n,0)^2\>
 -{\lambda_-\over 4} \min(x_n,0)^2 .$
Note that the nomenclature ``branch point'' is abusive in the
sense that it does {\em not necessarily} imply a bifurcation
of the
free boundary at that point (see Figure \ref{branchp}). Also there {\em are} one-phase
bifurcation points of the free boundary that are not included in
our class of branch points. Nevertheless it makes sense to speak
of branch points because {\em generically} a bifurcation occurs at
those points.
\\
In this paper we prove (cf. Theorem \ref{branch}) that {\em in two dimensions}
the free boundary is in a neighborhood
of each branch point the union of (at most) two $C^1$-graphs.
As application we obtain the following stability result:
If the free boundary contains no singular one-phase point
for certain boundary data $(B_0)$, then for boundary data $(B)$ close
to $(B_0)$ the free boundary consists of $C^1$-arcs
converging to those of $(B)$ (cf. Theorem \ref{stab}).
\\
In higher dimensions we derive an estimate for the
$(n-1)$-dimensional Hausdorff measure of the free boundary.
\\
Unfortunately the known techniques seem to be insufficient
to do a conclusive analysis at branch points. One reason
is that the density of the monotonicity formula
by H.W. Alt-L.A. Caffarelli-A. Friedman takes the value
$0$ at branch points.\\
The situation is complicated by the fact that the
limit manifold of all possible blow-ups at branch points
(including the case of varying centers) is not a
one-dimensional or even smooth manifold, but has a
more involved structure. Also the convergence to
blow-up limits is close to the branch-point {\em not
uniform!} Here we use an intersection-comparison approach
based on the Aleksandrov reflection
to show that -- although the flow with respect to the limit
manifold may not slow down when blowing up --
the free boundaries are still {\em uniformly graphs}
(see Proposition \ref{turn}).
The approach in Proposition \ref{turn} uses -- apart
from the reflection invariance -- very little information
about the underlying PDE and so yields a general
approach to the regularity of free boundaries in two
space dimensions provided that there is some information
on the blow-up limits.
\\
The Aleksandrov reflection has been recently used to prove
regularity in geometric parabolic PDE (\cite{chow1}, \cite{chow2}, \cite{chow3}).
In contrast to those results, where structural conditions for
the initial data are preserved under the flow, our
results are completely local.\\
{\bf Acknowledgment:} We wish to thank Nina Uraltseva for valuable
discussions and suggestions.
\section{Notation and Technical tools}
Throughout this article $\R^n$ will be equipped with the Euclidean
inner product $x\cdot y$ and the induced norm $\vert x \vert\> .$
$B_r(x)$ will denote the open $n$-dimensional ball of center
$x\> ,$ radius $r$ and volume $r^n\> \omega_n\> .$
When the center is not specified, it is assumed to be $0.$

We will use $\partial_e u = \nabla u\cdot e$ for the directional
derivative.

When considering a set $A\> ,$ $\chi_A$ shall stand for
the characteristic function of $A\> ,$
while
$\nu$ shall typically denote the outward
normal to a given boundary.

\label{exist}
Let $\lambda_+> 0$ and $\lambda_->0\> , \>
n \ge 2,$ let $\Omega$ be a bounded open subset of $\R^n$
with Lipschitz boundary and assume that $u_D \in W^{1,2}(\Omega)\> .$
From \cite{interf} we know then that there exists a ``solution'',
i.e. a function
$u\in W^{2,2}(\Omega)$
solving the strong equation $\Delta u = {\lambda_+ \over 2}
\> \chi_{\{u>0\}}\> - \>{\lambda_- \over 2}
\> \chi_{\{u<0\}}
$ a.e. in $\Omega$, and attaining the boundary data $u_D$ in
$L^2\> .$
The boundary condition may be replaced by other, more general
boundary conditions.

The tools at our disposition include two powerful monotonicity
formulae. One is the monotonicity formula introduced in
\cite{cpde} by one of the authors for a class of semilinear free
boundary problems (see also \cite{inv}). The second monotonicity
formula has been introduced by H.W. Alt-L.A. Caffarelli-A.
Friedman in \cite{mono}.
 What we are actually going to apply in section \ref{global} is a stronger
statement than the one in \cite{mono}.

For the sake of completeness let us state both monotonicity formulae here.
\begin{theorem}[Weiss's Monotonicity Formula]
\label{wmon}
Suppose that $B_\delta(x_0)\subset \Omega\> .$
Then for all $0<\rho<\sigma<\delta$
the function
\[ \Phi_{x_0}(r) := r^{-n-2} \int_{B_r(x_0)} \left(
{\vert \nabla u \vert}^2 \> +\> \lambda_+\max(u,0)
\> + \> \lambda_-\max(-u,0)\right)\]\[
- \; 2 \> r^{-n-3}\>  \int_{\partial B_r(x_0)}
u^2 \> d{\mathcal H}^{n-1}\; ,\]
defined in $(0,\delta)\> ,$ satisfies the monotonicity formula
\[ \Phi_{x_0}(\sigma)\> -\> \Phi_{x_0}(\rho) \; = \;
\int_\rho^\sigma r^{-n-2}\;
\int_{\partial B_r(x_0)} 2 \left(\nabla u \cdot \nu - 2 \>
{u \over r}\right)^2 \; d{\mathcal H}^{n-1} \> dr \; \ge 0 \; \; .\]
\end{theorem}
For a proof see \cite{cpde}.

In section \ref{global} we are going to need the following
stronger version of the Alt-Caffarelli-Friedman monotonicity
formula.

\begin{theorem}[Alt-Caffarelli-Friedman Monotonicity Formula]
\label{mon} Let $h_1$ and $h_2$ be continuous non-negative subharmonic
$W^{1,2}$-functions in $B_R(z)$ satisfying $h_1h_2=0$ in $B_R(z)$
as well as $h_1(z)=h_2(z)=0\> .$\\
Then for
\[ \Psi_z(r,h_1,h_2) := r^{-4} \int_{B_r(z)} {\vert \nabla h_1(x)\vert^2\over
{\vert x-z\vert^{n-2}}} \> dx \; \int_{B_r(z)} {\vert \nabla
h_2(x)\vert^2\over {\vert x-z\vert^{n-2}}} \> dx\; ,\]
and for $0<\rho<r<\sigma<R$, we have $\Psi_z(\rho)\leq
\Psi_z(\sigma)$. Moreover, if equality holds for some
$0<\rho<r<\sigma<R$
then one of the following is true:\\
(A) $h_1=0$ in $B_\sigma(z)$ or $h_2=0$ in $B_\sigma(z)$,\\
(B) for $i=1,2$, and $\rho<r<\sigma$,
 $\hbox{supp }(h_i) \cap \partial B_r(z)$ is a half-sphere and
 $h_i \Delta h_i=0$ in $B_\sigma(z)\setminus B_\rho(z)$ in the sense of measures.
\end{theorem}
For a proof of this version of monotonicity see \cite{SUW}. We
also refer to \cite{mono}, for the original proof.

It is noteworthy that
$$\Psi_z(r,(\partial_e u)^+, (\partial_e u)^-)=
\Psi_0(1,(\partial_e u_r)^+, (\partial_e u_r)^-) \text{ and } \Phi_z(r,u)=\Phi_0(1,u_r),
$$
where
$$u_r(x)=\frac{u(rx+z)}{r^2}.
$$

It is in fact possible
to apply Theorem \ref{mon} to the positive and
negative part of directional derivatives of $u\> :$
due to N. Uraltseva,
the functions $\max(\partial_e u,0)$ and $-\min(\partial_e
u,0)$ are subharmonic in $\Omega$ (see Lemma 2 in
\cite{uraltseva}).

A quadratic growth estimate near the set $\Omega\cap \{ u=0\} \cap
\{\nabla u=0\}$ had already been proved in \cite{interf} for more
general coefficients $\lambda_+$ and $\lambda_-\> ,$ but local
$W^{2,\infty}$- or $C^{1,1}$-regularity of the solution has been
shown for the first time in \cite{uraltseva}. See also
\cite{shah}. So we know that
\begin{equation}\label{C11}
u\in W^{2,\infty}_{\rm loc}(\Omega)\> .
\end{equation}

\begin{lemma}\label{homog}
Let $u$ be a solution of (\ref{obst}) in $B_1$ and
suppose that the origin is a free boundary point.
Then the following statements are equivalent:\\
1) Either $\nabla u(0)\ne 0,$ or $\lim_{r\to 0} \Psi_0(r,(\partial_e u)^+, (\partial_e u)^-)=0$
for each direction $e.$\\
2) Either $\nabla u(0)\ne 0,$ or each blow-up limit
$$u_0(x) = \lim_{m\to \infty} {u(r_m x)\over {r_m^2}}$$
is after rotation of the form
$$u_0(x) =
a_1 {\lambda_+\over 4} \max(x_1,0)^2\>  -a_2 {\lambda_-\over 4}
\min(x_1,0)^2$$
where $a_1,a_2\in \{0,1\}$ and $a_1+a_2\ne 0$.\\
3) Either $\nabla u(0)\ne 0,$ or at least one blow-up limit
$$u_0(x) = \lim_{m\to \infty} {u(r_m x)\over {r_m^2}}$$
is after rotation of the form
$$u_0(x) =
a_1 {\lambda_+\over 4} \max(x_1,0)^2\>  -a_2 {\lambda_-\over 4}
\min(x_1,0)^2$$
where $a_1,a_2\in \{0,1\}.$
\\
4) The origin is not a one-phase singular free boundary point,
i.e. no blow-up limit
$$u_0(x) = \lim_{m\to \infty} {u(r_m x)\over {r_m^2}}$$
is allowed to be a non-negative/non-positive homogeneous polynomial
of degree $2.$
\end{lemma}
\proof ``$1) \Rightarrow 2):$''
In the case $\nabla u(0)\ne 0,$ we obtain -- using
for example the lower semicontinuity of the weighted
$L^2$-norm $f\mapsto \int_{B_1} |x|^{2-n}f^2(x)\> dx$
with respect to weak convergence -- that
\[ 0 = \lim_{m\to\infty} \Psi_0(1,(\partial_e
u_{r_m})^+, (\partial_e u_{r_m})^-)
\ge \int_{B_1} {\vert \nabla (\partial_e
u_{0})^+(x)\vert^2\over
{\vert x\vert^{n-2}}} \> dx \; \int_{B_1} {\vert \nabla (\partial_e
u_{0})^-(x)\vert^2\over
{\vert x\vert^{n-2}}} \> dx\; .\]
Thus Theorem \ref{mon} (A) applies, and
we obtain that for each
direction $e,$ either $\partial_e u_0\ge 0$ in $\R^n$ or
$\partial_e u_0\le 0$ in $\R^n.$ It follows that after rotation,
$u_0$ is a function depending only
on the $x_1$ variable, and we obtain 2).\\
``$2) \Rightarrow 3)$'' is trivial.
``$3) \Rightarrow 1)$'' holds because the function in 2)
is one-dimensional and because the limit
$\lim_{r\to 0} \Psi_0(r,(\partial_e u)^+, (\partial_e u)^-)=0$ exists.\\
``$3) \Leftrightarrow 4):$''
From the monotonicity formula
\ref{wmon} (cf. \cite[Theorem 4.1]{cpde}) it follows that
in the case $\nabla u(0)= 0,$
$u_0$ is a $2$-homogeneous solution of the same equation.
These solutions have been characterized (cf. \cite[Theorem 4.3]{SUW},
and the only possibilities are
the solutions in 2) {\em and}
certain non-negative/non-positive homogeneous polynomials
of degree $2.$
\section{Classification of Global Solutions}
\label{global}
In what follows,
$I$ shall be an index set in a metric space.\\
We define the class
\begin{equation}\label{M^*}
\begin{split}
M^* &:= \{ u:B_1(0)\to \R \> :\\
&u(x_1,\dots,x_n) = \beta_1\left({\lambda_+\over 4} \max(x_1,0)^2\>
-{\lambda_-\over
4} \min(x_1-\tau,0)^2\right) + \beta_2 x_1,\\
&\text{where }\tau \in [-1,0], 0\leq \beta_1\leq a, 0\leq \beta_2
\leq b,
0<c\leq \beta_1+\beta_2,\\
&\text{ and } \beta_2\neq 0 \text{ implies } \tau=0\}.
\end{split}
\end{equation}
The class $M$ is then defined as
all rotated elements of $M^*$, i.e.
\begin{equation}\label{M}
M := \{ u:B_1(0)\to \R \> :\\
u= v \circ U \text{ where } U \text{ is a rotation, } v \in
M^*\}.
\end{equation}
Observe that singular one-phase solutions are excluded from
$M$.
\begin{theorem}\label{global-solutions}
Let $(u_\alpha)_{\alpha \in I}$ be a family of solutions of
(\ref{obst}) in $B_1$ that is bounded in $W^{2,\infty}(B_1),$
and suppose that $0\in \Omega\cap (\partial\{u_{\alpha_0}>0\}\cup
\partial\{u_{\alpha_0}<0\})$
for some ${\alpha_0}\in I$, and either
$\nabla u_{\alpha_0}(0)\ne 0$ or $\lim_{r\to 0} \Psi_0(r,(\partial_e
u_{\alpha_0})^+, (\partial_e u_{\alpha_0})^-)=0$ for each
direction $e;$ this means by Lemma
\ref{homog} that $0$ is
not a singular one-phase free boundary point. Define further $S_r$
by
$$
r^{n-1}S^2_r(y,u_\alpha)=\int_{\partial B_r(y)} u_\alpha^2,
$$
Then, if $u_\alpha \to u_{\alpha_0}$ in $L^1(B_1)$
as $\alpha \to \alpha_0,
\partial \{u_{\alpha}>0\}\ni  \ y \to 0$ and $r\to 0$, all possible limit
functions of the family
$$
\frac{u_{\alpha}(y + r \cdot)}{S_{r}(y,u_{\alpha})} ,
$$
belong to $M$ for some $a,b,c$ as above.
\end{theorem}

\proof
As the statement holds by the implicit function theorem in the case
$\nabla u_{\alpha_0}(0)\ne 0$, we may assume
$\nabla u_{\alpha_0}(0)=0$ and
$\lim_{r\to 0} \Psi_0(r,(\partial_e
u_{\alpha_0})^+, (\partial_e u_{\alpha_0})^-)=0$ for each
direction $e$.
Consider sequences $u_j:=u_{\alpha_j}\to u_{\alpha_0}$, $\partial \{u_j>0\}\ni  \ y_j \to 0, r_j\to 0$ and scaled functions
$$
v_j(x)=\frac{u_j(y_j + r_j x)}{S_{r_j}(y_j,u_j)}.
$$
A straightforward analysis of the limits of $v_j$ will yield the statement of our theorem.
First, setting
$$
T_j:=\frac{r_j^2 }{S_{r_j}(y_j,u_j)},
$$
we see that $T_j$ is uniformly bounded from above, due to the
non-degeneracy \cite[Lemma 3.7]{SUW}. Next,  by the bounds on the second derivatives,
$$
|D^2v_j(x)| = \frac{r_j^2}{S_{r_j}(y_j,u_j)}
|D^2u_j(y_j+r_jx)|\leq CT_j \leq C_0, \qquad x\in B_{1/(2r_j)},
$$
so  that the $W^{2,\infty}$-norm of $v_j$ is locally uniformly
bounded. Now as the free boundary has zero Lebesgue measure
\cite[Theorem 5.1]{interf} one can infer as in \cite{CKS}, General Remarks,  that
$v_j$ has a subsequence converging strongly
in $W_{\rm loc}^{2,p} (\R^n)$. Let $v$ be a
limit function.
The assumption $\lim_{r\to 0} \Psi_0(r,(\partial_e
u_{\alpha_0})^+, (\partial_e u_{\alpha_0})^-)=0$ implies now
by the monotonicity formula Theorem \ref{mon} that
for each $R\in (0,\infty)$ and $\delta>0,$
$$\delta \geq \Psi_{y_j}(Rr,(\partial_e u_j)^+,
(\partial_e u_j)^-)\geq \Psi_{y_j}(Rr_j,(\partial_e u_j)^+,
(\partial_e u_j)^-)$$ $$= \Psi_{0}(R,(\partial_e v_j)^+, (\partial_e
v_j)^-) \frac{S_{r_j}(y_j,u_j)}{r_j^2 }$$ 
if we choose first $r$ small and then $j$ sufficiently large.\\
Consequently $\Psi_{0}(R,(\partial_e v)^+, (\partial_e
v)^-)=0$ for every $R\in (0,\infty)$ and every direction $e.$ But
then
Theorem \ref{mon} (A) applies, and for each direction $e,$ either $\partial_e v\ge 0$ in $\R^n$
or $\partial_e v\le 0$ in $\R^n.$ In particular $v$ is one
dimensional.
As $\int_{\partial
B_1(0)}|v|^2 =\lim_{j\to\infty} \int_{\partial B_1(0)}|v_j|^2 =1,$
we obtain that $v\in M.$
\section{Uniform regularity of the free boundary close to branch points}
This chapter contains the main result of this paper.
\begin{theorem}\label{branch}
Let $n=2$,
let $(u_\alpha)_{\alpha \in I}$ be a family of solutions of
(\ref{obst}) in $B_1$ that is bounded in $W^{2,\infty}(B_1),$
and suppose that for some ${\alpha_0}\in I,$ a blow-up limit
$$\lim_{m\to\infty} {u_{\alpha_0}(r_m\cdot)\over {r_m^2}}$$
is contained in $M^*.$\\
Then, if $u_\alpha \to u_{\alpha_0}$ in $L^1(B_1)$ as
$\alpha \to \alpha_0,$ $B_{r_0}\cap \partial\{u_\alpha>0\}$ and
$B_{r_0}\cap
\partial\{u_\alpha<0\}$ are $C^1$-graphs uniformly in $\alpha \in
N_\kappa(\alpha_0)$ for some $r_0>0$ and $\kappa>0;$ here the
direction of every graph is the same,
and $N_\kappa(\alpha_0)$ is a given open neighborhood
of $\alpha_0$.
\end{theorem}

The crucial tool in the proof of the theorem is the following
proposition which uses an Aleksandrov reflection approach.

\begin{proposition}\label{turn}
Let $n=2$,
let $(u_\alpha)_{\alpha \in I}$ be a family of solutions of
(\ref{obst}) in $B_1$ that is bounded in $W^{2,\infty}(B_1),$
and suppose that
for some ${\alpha_0}\in I,$ a blow-up limit
$$\lim_{m\to\infty} {u_{\alpha_0}(r_m\cdot)\over {r_m^2}}$$
is contained in $M^*.$\\
Then, given $\epsilon\in (0,1/8)$
there exist positive $\kappa,\delta$ and $\rho$
such that for $\alpha \in N_\kappa(\alpha_0), y\in B_\delta \cap \partial\{u_\alpha>0\}$ and
$r\in (0,\rho),$ the scaled function
\begin{equation}\label{scale-at-y}
u_r(x)=\frac{u_\alpha(rx +y)}{S_r(y,u_\alpha)}
\end{equation}
satisfies
$$
\dist(u_r,M^*) = \inf_{v\in M^*}\sup_{B_1(0)} \left| v(x) - u_r(x)
\right| < \epsilon.
$$
\end{proposition}

The result implies
that we have uniform control of the rotation of the free
boundaries. In particular, this implies uniform cone-flatness
of the free boundaries.

\proof First, by Theorem \ref{global-solutions},
for any $\tilde \epsilon>0$ there are positive $\tilde \kappa,\tilde\delta$
and $\tilde\rho$ such that
$$ \dist(u_r,M) < \tilde \epsilon \text{ for }
\alpha \in N_{\tilde \kappa}(\alpha_0), y \in \partial\{ u_\alpha>0\} \cap
B_{\tilde\delta}\text{ and } r\in (0,\tilde \rho).
$$
Now if the statement of the theorem does not hold, then there are
positive $r_0$ and $r_1$ as well as two counterclockwise rotations
$U_{\theta_0}$ and $U_{\theta_1}$ of non-negative
angle
$\theta_0$ and $\theta_1,$ respectively, satisfying
$\vert\theta_0-\theta_1\vert \ge c_1\epsilon > 0$ and
$$ \dist(u_{r_0}\circ U_{\theta_0}, M^*) \leq \tilde \epsilon
\qquad \text{ as well as } \qquad \dist(u_{r_1}\circ U_{\theta_1},
M^*) \leq \tilde\epsilon;$$
here $c_1$ is a constant depending on $(a,b,\lambda_+,\lambda_-)$.\\
Let now $M^{*,\theta} := \{ v:B_1(0)\to \R : v\circ U_\theta \in
M^*\}$ 
and observe that while we do not know at this stage
whether $r \mapsto u_r$ is uniformly continuous on
$(0,\tilde \rho)$, we do know that
$t \mapsto u_{\exp(-t)}$ is uniformly continuous on 
$(t_0,+\infty)$. 
As each continuous connection of $M^{*,\theta_0}$ and
$M^{*,\theta_1}$ in $M$ must either contain for each $\theta\in
[\theta_0,\theta_1]$ an element of $M^{*,\theta},$ or contain for
each $\theta\in [-\pi,\pi)\setminus (\theta_0,\theta_1)$ an
element of $M^{*,\theta},$ we obtain for small $\tilde \epsilon$
--- depending on $(\epsilon,a,b,c,\lambda_+,\lambda_-)$
--- also $c_2 \in [c_1/4,3c_1/4]$ and
$0<r_2<r_3<1$ as well as two rotations
$U_{\theta_2}$ and $U_{\theta_3}$ satisfying $\vert\theta_3-\theta_2\vert = {c_2\epsilon}$
such that
$$ \dist(u_{r_2}\circ U_{\theta_2}, M^*) \leq \tilde \epsilon
\qquad \text{ and } \qquad  \dist(u_{r_3}\circ U_{\theta_3}, M^*)
\leq \tilde \epsilon.$$ We may assume that $\theta_3-\theta_2>0;$
if this is not the case, we apply the following part of the proof
to $u_\alpha(x_1,2y_2-x_2)$ instead of $u_\alpha(x_1,x_2).$ Now
set
$$\omega = {c_2\epsilon\over 2}, \qquad U =
U_{{\theta_2+\theta_3\over 2}}.$$
\begin{figure}
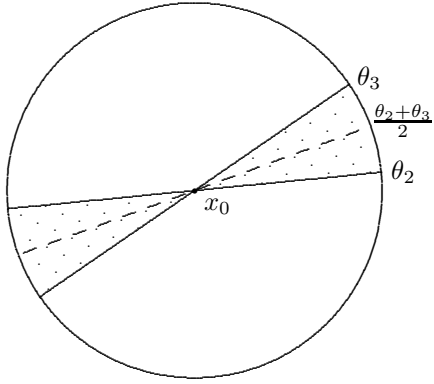

\[
\beginpicture
\setcoordinatesystem units <2.5cm,2.5cm> \setplotarea x from -1 to
1, y from 0 to 2 \startrotation by 0.937 0.3420 about 0 1 \put
{\makebox(0,0)[t]{$x_0$}} [Bl] at 0.1 0.92 \put {\circle*{2}} [Bl]
at 0 1 \linethickness1pt \put {\makebox(0,0)[t]{$\theta_3$}} [Bl]
at 1.1 1.30 \put {\makebox(0,0)[t]{$\theta_2$}} [Bl] at 1.1 0.77
\put {\makebox(0,0)[t]{${\theta_2+\theta_3\over 2}$}} [Bl] at 1.20
1.05 \linethickness2pt \plot -0.969 1.25 0.969 0.75 / \plot -0.969
0.75 0.969 1.25 / \setdashes \plot -0.969 1 0.969 1 / \setsolid
\linethickness2mm \circulararc 360 degrees from 1 1 center at 0 1
\setlinear \vshade -0.969 0.75 1.25 0 1 1 / \vshade 0 1 1 0.969
0.75 1.25 /
\endpicture
\]
\\[.20cm]
\begin{center}
\caption{Turning free boundary}\label{rotations}
\end{center}
\end{figure}
Moreover, let
 $$
 \phi (r,\theta):=\frac{u_\alpha(y + rU(\cos \theta,\sin
 \theta))}{S_r(y,u_\alpha)}\; .
 $$
For each $0<r<1/2$, the function $\phi(r,\cdot
 )$ defines a function on the unit circle $[-\pi,\pi).$
The following part is
inspired by applications of the Aleksandrov reflection (see for example
\cite{chen}, \cite{matanoveron}, \cite{veronveron}).
There are however important differences: while the authors in
\cite{chen}, \cite{matanoveron}, \cite{veronveron} exclude
{\em repetitive} behavior as $r\to 0$,
for our application it is necessary to derive
a contradiction from {\em just one turn of angle}
$|\theta_0-\theta_1|$. Moreover, our class $M$
is not a one-dimensional or even a smooth manifold.\\
We consider
$$
 \xi (r,\theta):=\phi(r,\theta) - \phi(r,-\theta)
 $$
and observe that $\xi(r,0)=\xi(r,\pi)=0$. In what follows we will
prove that $\xi(r_3,\theta) \ge 0$ for $\theta \in [0,\pi]$ and
${\partial\xi\over {\partial\theta}}(r_2,0) \; < \; 0$ provided
that $\tilde\epsilon$ has been chosen small enough (depending on
$(\epsilon,a,b,c,\lambda_+,\lambda_-,\sup_{\alpha\in I}\sup_{B_1(0)} \vert u_\alpha
\vert)$). By the comparison principle (applied to
$S_r(y,u_\alpha)\phi(r,\theta)$ and $S_r(y,u_\alpha)\phi(r,-\theta)$ in the
two-dimensional domain $[0,r_3)\times (0,\pi)$ with respect to the
original coordinates $x_1$ and $x_2$) this yields a contradiction.
\\
Let us prove
${\partial\xi\over {\partial\theta}}(r_3,0)
\; > \; 0$ as well as
$\xi(r_3,\theta) \ge 0$
for $\theta \in [0,\pi].$
The partial derivative estimate at $r_2$ is obtained in the same way.
\begin{figure}
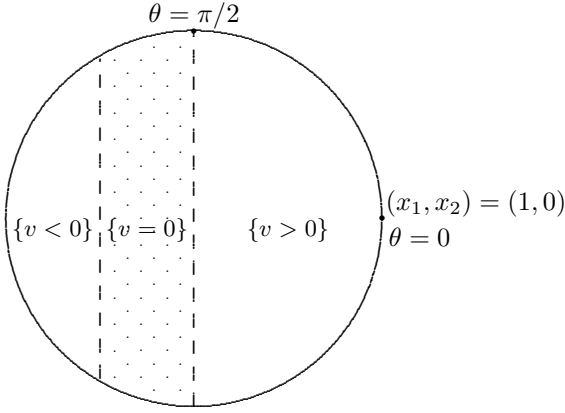

\[
\beginpicture
\setcoordinatesystem units <2.5cm,2.5cm>
\setplotarea x from -1 to 1, y from -1 to 1.5
\put {\makebox(0,0)[t]{{\small $\{ v>0\}$}}} [Bl] at 0.5 0
\put {\makebox(0,0)[t]{{\small $\{ v=0\}$}}} [Bl] at -0.25 0
\put {\makebox(0,0)[t]{{\small $\{ v<0\}$}}} [Bl] at -0.75 0
\put {\circle*{2}} [Bl] at 1 0
\put {\circle*{2}} [Bl] at 0 1
\linethickness1pt
\put {\makebox(0,0)[t]{$(x_1,x_2)=(1,0)$}} [Bl] at 1.5 0.15
\put {\makebox(0,0)[t]{$\theta=0$}} [Bl] at 1.2 -0.05
\put {\makebox(0,0)[t]{$\theta=\pi/2$}} [Bl] at 0 1.15
\linethickness2pt
\setdashes
\plot 0 -1 0 1 /
\plot -.5 -0.85 -.5 0.85 /
\setsolid
\linethickness2mm \circulararc 360 degrees from 0 1 center at 0 0
\setlinear
\vshade -.5 -0.85 0.85 0 -1 1 /
\endpicture
\]
\\[.20cm]
\begin{center}
\caption{Example of $v$}\label{vexamp}
\end{center}
\end{figure}
Take $v\in M^*$ such that
$$\sup_{B_1(0)} \vert v - u_{r_3}\circ U_{\theta_3}\vert
= \dist(u_{r_3}\circ U_{\theta_3},M^*) \le \tilde\epsilon$$ (note
that we do not need the axiom of choice to do so) and define
$$\phi_0(\theta) := v(\cos \theta,\sin \theta),$$
$$\sigma := \sup\{ \theta \in (0,\pi):\phi_0(\theta)=0\}$$
$$\text{and } \xi^\gamma_0(\theta):=\phi_0(\gamma+\theta) - \phi_0(\gamma-\theta).$$
\begin{figure}
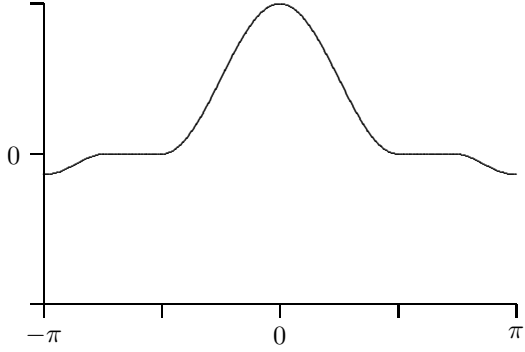

\[
\beginpicture
\setcoordinatesystem units <1cm,1cm>
\setplotarea x from -3.14159 to 3.14159, y from -2 to 2
\axis bottom ticks quantity 5 /
\axis left ticks quantity 3 /
\put {\makebox(0,0)[t]{$-\pi$}} [Bl] at -3.14159 -2.3
\put {\makebox(0,0)[t]{$\pi$}} [Bl] at 3.14159 -2.3
\put {\makebox(0,0)[t]{$0$}} [Bl] at 0 -2.3
\put {\makebox(0,0)[t]{$0$}} [Bl] at -3.54159 0.1
\linethickness10pt
\setquadratic
\plot -3.14159 -0.27596
-3.01593 -0.261798
-2.89027 -0.22175
-2.7646 -0.162965
-2.63894 -0.0967684
-2.51327 -0.0378408
-2.38761 -0.00312249
-2.26195 0
-2.13628 0
-2.01062 0
-1.88496 0
-1.75929 0
-1.63363 0
-1.50796 0.00788529
-1.3823 0.0702235
-1.25664 0.190983
-1.13097 0.362576
-1.00531 0.574221
-0.879646 0.812619
-0.753982 1.06279
-0.628319 1.30902
-0.502655 1.53583
-0.376991 1.72897
-0.251327 1.87631
-0.125664 1.96858
0 2
0.125664 1.96858
0.251327 1.87631
0.376991 1.72897
0.502655 1.53583
0.628319 1.30902
0.753982 1.06279
0.879646 0.812619
1.00531 0.574221
1.13097 0.362576
1.25664 0.190983
1.3823 0.0702235
1.50796 0.00788529
1.63363 0
1.75929 0
1.88496 0
2.01062 0
2.13628 0
2.26195 0
2.38761 -0.00312249
2.51327 -0.0378408
2.63894 -0.0967684
2.7646 -0.162965
2.89027 -0.22175
3.01593 -0.261798
3.14159 -0.27596 /
\endpicture
\]
\caption{Example of $\phi_0$}\label{halfplane}
\end{figure}
\begin{figure}
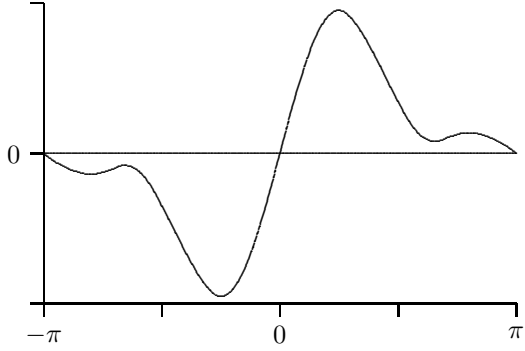

\[
\beginpicture
\setcoordinatesystem units <1cm,1cm>
\setplotarea x from -3.14159 to 3.14159, y from -2 to 2
\axis bottom ticks quantity 5 /
\axis left ticks quantity 3 /
\put {\makebox(0,0)[t]{$-\pi$}} [Bl] at -3.14159 -2.3
\put {\makebox(0,0)[t]{$\pi$}} [Bl] at 3.14159 -2.3
\put {\makebox(0,0)[t]{$0$}} [Bl] at 0 -2.3
\put {\makebox(0,0)[t]{$0$}} [Bl] at -3.54159 0.1
\linethickness2pt
\plot -3.14159 0.0 3.14159 0.0 /
\linethickness10pt
\setquadratic
\plot -3.14159 0
-3.01593 -0.0936459
-2.89027 -0.162964
-2.7646 -0.22175
-2.63894 -0.261798
-2.51327 -0.27596
-2.38761 -0.261798
-2.26195 -0.22175
-2.13628 -0.17085
-2.01062 -0.166992
-1.88496 -0.228824
-1.75929 -0.365698
-1.63363 -0.574221
-1.50796 -0.812619
-1.3823 -1.06279
-1.25664 -1.30902
-1.13097 -1.53583
-1.00531 -1.72897
-0.879646 -1.86842
-0.753982 -1.89836
-0.628319 -1.80902
-0.502655 -1.60601
-0.376991 -1.30209
-0.251327 -0.91635
-0.125664 -0.473036
0 0
0.125664 0.473036
0.251327 0.91635
0.376991 1.30209
0.502655 1.60601
0.628319 1.80902
0.753982 1.89836
0.879646 1.86842
1.00531 1.72897
1.13097 1.53583
1.25664 1.30902
1.3823 1.06279
1.50796 0.812619
1.63363 0.574221
1.75929 0.365698
1.88496 0.228824
2.01062 0.166992
2.13628 0.17085
2.26195 0.22175
2.38761 0.261798
2.51327 0.27596
2.63894 0.261798
2.7646 0.22175
2.89027 0.162964
3.01593 0.0936459
3.14159 0 /
\endpicture
\]
\caption{Example of $\xi^\gamma_0$}\label{xi0}
\end{figure}
Observe that we may assume $\sigma = \pi$ or
$\sigma\le 3\pi/4.$ If this is not the case we change $r_3$
to $r_3/2$ where we still have
flatness in the same direction.
\\
Since $\phi_0$ is an even $2\pi$-periodic function which is
decreasing on $(0,\pi),$ $\xi^\gamma_0(\theta) \ge 0$ for $-\pi\le
\gamma \le 0, 0\le \theta \le \pi.$ Note that in the case of
$\phi_0$ being strictly decreasing in $(0,\pi),$ we also obtain
$\xi^\gamma_0(\theta) > 0$ for $-\pi < \gamma < 0, 0 < \theta < \pi.$
As $u_{r_3}\circ U_{\theta_3}$ is close to $v$ (and thus
$\xi(r_3,\cdot)$ close to $\xi^{-\omega}_0$)
we expect the same to hold for
$\xi(r_3,\cdot).$
In order to prove this rigorously we proceed as follows:\\
1) ${\partial\over
{\partial\theta}}\xi^\gamma_0(0)=2\phi'_0(\gamma) \ge
c_3=c_3(\gamma,c,\lambda_+,\lambda_-)>0$ for $\gamma \in
(-\pi/2,0),$ and for $\sigma\ne \pi,$ ${\partial\over
{\partial\theta}}\xi^\gamma_0(\pi)=2\phi'_0(\gamma+\pi) \le
-c_3<0$ for $\gamma \in (-\pi+\sigma,0).$ Consequently, for small
$\tilde \epsilon$ (depending on
$(\epsilon,a,b,c,\lambda_+,\lambda_-,\sup_{\alpha\in
I}\sup_{B_1(0)} \vert u_\alpha \vert)$), ${\partial\over
{\partial\theta}}\xi(r_3,0) \ge c_3/2>0$ and, in the case
$\sigma\ne \pi,$ ${\partial\over {\partial\theta}}\xi(r_3,\pi) \le
-c_3/2<0.$ It follows that there is $c_4 =
c_4(\epsilon,a,b,c,\lambda_+,\lambda_-,\sup_{\alpha\in
I}\sup_{B_1(0)} \vert u_\alpha \vert)$ such that
$\xi(r_3,\cdot)>0$ in $(0,c_4)$ and, in the case $\sigma\ne \pi,$
$\xi(r_3,\cdot)>0$ in $(\pi-c_4,\pi).$\\
2) Next, since $\phi_0(\gamma+\theta)\ge 0$ for $\theta\in [-\sigma-\gamma,\sigma-\gamma]$
and $\phi_0(\gamma-\theta)\le 0$ for $\theta\in [\pi/2+\gamma,\pi+\gamma],$
making use of the non-degeneracy \cite[Lemma 3.7]{SUW},
we see that for small $\tilde \epsilon$
(depending on $(\epsilon,a,b,c,\lambda_+,\lambda_-,\sup_{\alpha\in I}\sup_{B_1(0)} \vert u_\alpha \vert)$),
$$\phi(r_3,\theta) \ge 0 \text{ for }0\le \theta \le \sigma+\omega-\omega/2$$
and
$$\phi(r_3,-\theta) \le 0 \text{ for }\pi-3\omega/2\ge \theta \ge \pi/2-\omega+\omega/2.$$
Consequently
$$\xi(r_3,\theta) \ge 0 \text{ for }\pi/2-\omega/2\le \theta \le \sigma+\omega/2$$
and $\tilde \epsilon$ as above. Observe that $\sigma=\pi$ and
$\pi/2\leq \sigma\leq 3\pi/4$ are both allowed here.\\
3) Last, in $[c_4, \pi/2-\omega/4]\cup [\sigma+\omega/4,\pi-c_4],$
we obtain by the assumed range for $\sigma$ that
$$\xi_0^{-\omega}(\theta) = \phi_0(\theta-\omega)-\phi_0(\theta+\omega)
\geq c_5 = c_5(\epsilon,c,\lambda_+,\lambda_-)>0,$$
so that $\xi(r_3,\cdot) \ge c_5/2 >0$ in
$[c_4, \pi/2-\omega/4]\cup [\sigma+\omega/4,\pi-c_4]$ for small
$\tilde \epsilon.$

Combining 1)-3) we obtain the desired estimate, i.e.
${\partial\xi\over {\partial\theta}}(r_3,0)
\; > \; 0$ as well as
$\xi(r_3,\theta) \ge 0$
for $\theta \in [0,\pi].$
\\[.5cm]
{\em Proof of Theorem \ref{branch}:} By Proposition \ref{turn} we
know that $g^+_\alpha,g^-_\alpha$ defined by
$$g^+_\alpha(x_2) = \sup\{ x_1 \> : \> (x_1,x_2)\in B_\delta\cap \{u_\alpha=0\}\}$$
$$\text{and }g^-_\alpha(x_2) = \inf\{ x_1 \> : \> (x_1,x_2)\in B_\delta\cap \{u_\alpha=0\}\}$$
are
bounded in $C^1([-\delta/2,\delta/2]).$\\
We maintain that $u_\alpha=0$ in $B_{\tilde \delta}\cap \{
g^-_\alpha < x_n < g^+_\alpha\}$ for $\alpha \in
N_{\tilde\kappa}(\alpha_0).$ Suppose this is not true: then,
replacing if necessary $u$ by $-u$ and exchanging $\lambda_+$ and
$\lambda_-,$ there are $y,z\in B_{\tilde \delta}\cap
\partial\{ u_\alpha>0\}$
such that $y_2=z_2$ and $u_\alpha>0$ on
the straight line
segment between $y$ and $z.$ But then, setting
$r=2\vert y-z\vert,$ we obtain that
$$u_r(x) = \frac{u_\alpha(rx +y)}{S_r(y,u_\alpha)}$$
does not satisfy
$\dist(u_r,M^*) \le \tilde \epsilon,$ a contradiction to
Proposition \ref{turn} provided that $\tilde \kappa$ and
$\tilde \delta$ have been chosen small enough.

\section{Stability of the free boundary}
\begin{theorem}\label{stab}
Let $\Omega\subset\R^2$ be a bounded Lipschitz domain
and assume that for given Dirichlet data $u_D\in W^{1,2}(\Omega)$
the free boundary does not contain any one-phase
singular free boundary point (cf. Lemma \ref{homog}).\\
Then for $K\subset\subset \Omega$ and $\tilde u_D\in
W^{1,2}(\Omega)$ satisfying $\sup_{\partial\Omega} \vert
u_D-\tilde u_D\vert <\delta_K,$
there is $\omega>0$ such that
the free boundary is for every $y\in K$ in $B_\omega(y)$
the union of (at most) two $C^1$-graphs
which approach those of the
solution with respect to boundary data $u_D$ as
$\sup_{\partial\Omega} \vert u_D-\tilde u_D\vert\to 0.$
\end{theorem}
\proof
Let $u$ and $\tilde u$ be the solutions with respect to $u_D$ and
$\tilde u_D$, respectively. By the comparison principle,
$\sup_{\Omega} \vert u-\tilde u\vert\to 0$ as
$\sup_{\partial\Omega} \vert u_D-\tilde u_D\vert\to 0.$
Consequently, $\tilde u\to u$ in $C^{1,\beta}_{\rm loc}(\Omega)$ as
$\sup_{\partial\Omega} \vert u_D-\tilde u_D\vert\to 0.$
But then Theorem \ref{branch} applies, and the free boundary
of $\tilde u$ is in $B_\omega(y)$ the union of two
$C^1$-graphs which are bounded in $C^1.$ More precisely,
fixing $z\in\Omega\cap (\partial\{ u>0\}\cup \partial\{ u<0\})$
and translating and rotating once, we obtain $r_0>0$ such that
$\partial\{ \tilde u>0\}\cup \partial\{ \tilde u<0\}$
is for $\sup_{\partial\Omega} \vert u_D-\tilde u_D\vert<\delta_K$ in $B_{r_0}$
the union of the graphs of the $C^1$-functions $\tilde g^+$ and
$\tilde g^-$ in the direction of $e_2;$ moreover,
the $C^1$-norms of $\tilde g^+$ and
$\tilde g^-$ are bounded as $\tilde u_D\to u_D.$
Suppose now that
$$\sup_{[-r_0/2,r_0/2]} \vert \tilde g^+-\tilde g^-\vert \ge c_1 >0$$
for some sequence $\tilde u_D\to u_D.$ Then the fact that $u$ and
$\tilde u$ are near free boundary points after rescaling close to
$M$ (Theorem \ref{global-solutions}), implies that
$$\sup_{B_{1/2}} \vert \tilde u-u\vert \ge c_2 >0$$
for the same sequence, and we obtain a contradiction.
\section{Finite $n-1$-dimensional Hausdorff measure of the free boundary}

In this section we assume that $n\ge 2.$

We first show that the free
boundary has finite perimeter, which can be done as in \cite{brezis}:
Set
$$\beta (u):=\lambda_+ \chi_{\{u>0\}} -\lambda_- \chi_{\{u<0\}}$$
and define
$$\psi_\epsilon (t)=\left\{\begin{array}{ll}
1&\text{for } t>\epsilon,\\
-1&\text{for } t<-\epsilon, \text{ and}\\
t/\epsilon&\text{when }-\epsilon \leq t\leq
\epsilon.\end{array}\right.$$

Now, if $\eta$ is a cut-off function, we obtain, differentiating the
equation $\Delta u=\beta (u)$, multiplying by $\psi_\epsilon
(\partial_i u) \eta$ and integrating over $\Omega$, that
$$
\int_\Omega \beta' (u) \partial_i u\psi_\epsilon (\partial_i u)
\eta =\int_\Omega  \partial_i \Delta u\psi_\epsilon (\partial_i u)
\eta = -\int_\Omega \psi_\epsilon' |\nabla \partial_i u|^2\eta
-\int_\Omega \psi_\epsilon(\partial_i u) \partial_i \nabla u \cdot \nabla \eta .
$$
The first integral on the right-hand side of the equality being
non-positive and the second one bounded implies, letting
$\epsilon$ tend to zero, that
$$ \int_\Omega |\nabla \beta (u)|\eta
\leq C_1\; .
$$
Here we used the fact that $\psi_\epsilon$ converges to the sign
function as $\epsilon\to 0,$ and that $\beta' \geq 0$.
The above calculation can be made rigorous regularizing the
equation by $\Delta u_\delta = \beta_\delta(u_\delta)$
where $\beta_\delta$ is a smooth increasing function tending
to $\beta$ as $\delta\to 0;$ we let first $\epsilon$ and then
$\delta$ go to $0.$\\

Using in the above regularization the assumption $\min (\lambda_+,\lambda_-)>0$
as well as the lower semicontinuity of the $BV$-norm,
we obtain that
the sets $\{u>0\}$ and $\{u<0\}$ are locally in $\Omega$ sets of finite
perimeter. Since the set $\{u=0\}\cap \{\nabla u \neq 0\}$ is
locally in $\Omega$ a $C^{1,1}$-surface, the finite perimeter
estimate tells us that
$$
{\mathcal H}^{n-1}\left( \{u=0\}\cap \{\nabla u \neq 0\} \cap
K \right) <+\infty\; \hbox{ for each } K\subset\subset\Omega\; .
$$

Note that the above estimate implies also that
\begin{equation}\label{bound}
\int_\Omega |\nabla \Delta u|\eta\leq C_2 \int_\Omega \vert\nabla \eta\vert\; .
\end{equation}
This estimate in turn can be used to prove as in \cite{caff} that
$(\partial\{ u>0\}\cup\partial\{ u<0\})\cap \{\nabla u=0\}$ has locally in $\Omega$ finite
$n-1$-dimensional Hausdorff measure: for $\psi_\epsilon$ and $\eta$ as above,

\begin{equation}\label{identity}
\int_\Omega \eta \left(\nabla \psi_\epsilon(\partial_i u)\cdot \nabla {\partial_i u} +
\psi_\epsilon(\partial_i u) \Delta {\partial_i u}\right)= -\int_\Omega
\psi_\epsilon(\partial_i u)
\nabla \eta \cdot\nabla{\partial_i u}\; .
\end{equation}
Using estimate (\ref{bound}),
we deduce that
$$
\int_{\{0<|\partial_i u|<\epsilon\}\cap \Omega} \eta |\nabla {\partial_i u}|^2
\leq C_3 \epsilon \left( \int_{\Omega} \eta |\Delta {\partial_i u}| + \int_{\Omega} |D^2 u| \> |\nabla \eta|\right)\leq C_4 \epsilon \int_{\Omega} |\nabla \eta|.
$$
Take now -- using Vitali's covering theorem -- for each $\epsilon>0$ a covering $\bigcup_{j=1}^m B_\epsilon(x_j)$
of $(\partial\{ u>0\}\cup\partial\{ u<0\})\cap \{\nabla u=0\}\cap \{ \eta > 1\}$ such that
$x_j\in (\partial\{ u>0\}\cup\partial\{ u<0\})\cap \{\nabla u=0\}$ and
$\sum_{j=1}^m \chi_{B_\epsilon(x_j)}(y) \le C_5$ for
all $y\in \Omega;$ here $C_5$ depends only on the dimension $n.$
From the local $C^{1,1}$-regularity (\ref{C11}) and the non-degeneracy \cite[Lemma 3.7]{SUW}
we conclude as in \cite{caff} that
$$
{\mathcal L}^n(\{0<|\nabla u| <
\epsilon\}\cap B_\epsilon(x_j))\geq c_6 \epsilon^n
$$
where $c_6$ does not depend on
$\epsilon$ or $j.$ It follows that
$$
\sum_{j=1}^m \epsilon^{n-1}
\leq {1\over {c_6^n}} {1\over \epsilon} \sum_{j=1}^m
\left|\{0<|\nabla u| <\epsilon\} \cap B_\epsilon(x_j)\right| \leq
C_7 {1\over \epsilon} \sum_{j=1}^m
 \int_{\{0<|\nabla u| <\epsilon\} \cap B_\epsilon(x_j)} |\Delta u|^2 $$
 $$\leq C_8
  {1\over \epsilon} \sum_{j=1}^m\sum_{i=1}^n
 \int_{\{0<|\partial_i u| <\epsilon\} \cap B_\epsilon(x_j)} |\partial_{ii} u|^2
\leq C_9 {1\over \epsilon} \sum_{i=1}^n
\int_{\{0<|\partial_i u| <\epsilon\}} \eta |\partial_{ii} u|^2
\leq C_{10}
$$
where $C_{10}$ does not depend on $\epsilon.$
\\
We obtain:
\begin{theorem}\label{hausdorff}
Let $u$ be a solution of (\ref{obst}) in $\Omega.$ Then $\partial\{ u>0\}\cup
\partial\{ u<0\}$ is locally in $\Omega$ a set of finite
$n-1$-dimensional Hausdorff measure.
\end{theorem}

\end{document}